\documentclass[12pt]{amsart}
\usepackage{amsfonts}
\usepackage{amsmath,amssymb,amsthm}

\newtheorem{thm}{Theorem}[section]

\newtheorem{prop}[thm]{Proposition}

\newcommand{\R}{\mathbb R}

\newcommand{\E}{\mathbb E}

\newcommand{\X}{\mathfrak X}

\renewcommand{\span}{\mathrm{span}}

\newcommand{\ds}{\displaystyle}

\title{\bf Minimal surfaces in $S^3$ foliated by circles}
\author{N. Kutev \and V. Milousheva}
\address{Bulgarian Academy of Sciences, Institute of Mathematics and Informatics,
Acad. G. Bonchev Str. bl. 8, 1113 Sofia, Bulgaria} \email{}

\address{Bulgarian Academy of Sciences, Institute of Mathematics and Informatics,
Acad. G. Bonchev Str. bl. 8, 1113, Sofia, Bulgaria;
 "L. Karavelov" Civil Engineering
Higher School, 175 Suhodolska Str., Sofia, Bulgaria} \email{vmil@math.bas.bg}

\date{}



\begin{document}

\maketitle \thispagestyle{empty}

\begin{abstract}
We deal with minimal surfaces in the unit sphere $S^3$, which are
one-parameter families of circles. Minimal surfaces in $\R^3$
foliated by circles were first investigated by Riemann, and a
hundred years later Lawson constructed examples of such surfaces
in $S^3$. We prove that in $S^3$ there are only two types of
minimal surfaces foliated by circles, crossing the principal lines
at a constant angle. The first type surfaces are foliated by great
circles, which are bisectrices of the principal lines, and we show
that these minimal surfaces are the well-known examples of Lawson.
The second type surfaces, which are new in the literature, are
families of small circles, and the circles are principal lines. We
give a constructive formula for these surfaces. An application to
the theory of minimal foliated semi-symmetric hypersurfaces in
$\R^4$ is given.
\end{abstract}

\section{Introduction} \label{S:Intr}

In the present paper we deal with minimal surfaces in the unit
sphere $S^3$ in the four-dimensional Euclidean space $\R^4$,
equipped with the standard Euclidean metric $\langle ., .
\rangle$. A surface $M^2$ in $S^3$ is given by a unit
vector-valued function $l(u,v)$ in  $\R^4$, defined in a domain
${\mathcal D} \subset \R^2$, i.e.
$$l(u,v) = \left(l^1(u,v), l^2(u,v), l^3(u,v), l^4(u,v)\right), \,\, (u,v) \in {\mathcal D},$$
where $\langle l(u,v), l(u,v) \rangle = 1, \, (u,v) \in {\mathcal
D}$. Since our considerations are local, we assume that the
parameters $(u,v)$ are isothermal (conformal) ones, which means
that
 $\langle l_u, l_u \rangle = \langle l_v, l_v \rangle; \;
\langle l_u, l_v \rangle = 0$.

The minimal surfaces in $S^3$ are determined by the solutions $l =
l(u,v)$ of the following system of
 partial differential equations
$$\Delta l + |\nabla l|^2 \,l = 0; \leqno(1.1)$$
$$\langle l_u, l_u \rangle = \langle l_v, l_v \rangle; \qquad
\langle l_u, l_v \rangle = 0; \qquad \langle l, l \rangle = 1,$$
where $\nabla$ and $\Delta$ denote the gradient and  the laplacian
operators, respectively, computed with respect to the Euclidean
metric in $\R^4$.

System (1.1) is  the Euler - Lagrange system of harmonic maps,
which has been intensively investigated in the last decades by
variational methods (see for example \cite{G-M-S}, \cite{Hild},
\cite{Jost-1}, \cite{Jost-2}, \cite{Struwe} and the references
there).

Our aim is to find the minimal surfaces in $S^3$, that satisfy a
certain geometric property: locally they are one-parameter
families of circles. The variational methods cannot be applied for
studying the geometric structure of the minimal surfaces, that is
why we use a different method, which is based rather on the
differential geometry of surfaces in $\R^4$ than on PDE methods.

It is well known that the only minimal rotational surface in
$\R^3$ is the catenoid, which is a surface fibred by circles in
parallel planes. Other minimal surfaces in $\R^3$, which are
foliated by circles in parallel planes, were discovered by Riemann
\cite{Riemann}. They are usually referred in the literature as
\textit{Riemann examples}. Enneper \cite{Enneper} proved that
catenoids and Riemann examples are the only minimal surfaces in
$\R^3$ foliated by circles. A surface in $\R^3$, which is
determined by a smooth one-parameter family of circles is also
called a cyclic surface. Cyclic surfaces of constant mean
curvature and cyclic surfaces of constant Gauss curvature in
$\R^3$ are described in \cite{Nitsche} and \cite{Lopez}.

Our idea to find the minimal surfaces in $S^3$, which are
one-parameter family of circles, is motivated by what happens for
cyclic minimal surfaces in $\R^3$.

A well known example of a minimal surface in $S^3$ is the Clifford
torus (the standard flat torus), which consists of two orthogonal
families of circles. Generalizations of the Clifford torus are the
so called Lawson tori \cite{Lawson}. They have two orthogonal
families of parametric lines, one of them consists of circles, and
the other one consists of curves with constant Frenet curvatures
in $\R^4$.

The circles on the Lawson torus cross the principal lines at an
angle $\ds{\frac{\pi}{4}}$. In this paper we find all minimal
surfaces in $S^3$, which are one-parameter families of circles,
crossing the principal lines at a constant angle. We call these
surfaces \textit{generalized tori}. In Theorem \ref{T:generalized
torus} we prove that there are only two types of generalized tori
in $S^3$: the first one is characterized by the condition that the
circles are bisectrices of the  principal lines (we call these
surfaces  \textit{generalized tori of first type}), and the second
one is characterized by the condition that the circles are
principal lines (we call them \textit{generalized tori of second
type}). We show that all generalized tori of first type are Lawson
tori (Theorem \ref{T:First type}). In Theorem \ref{T:Second type}
we give a constructive formula for the generalized tori of second
type.

In Section \ref{S:Application} we point out the relation between
the theory of minimal surfaces in $S^3$ and the theory of minimal
foliated semi-symmetric hypersurfaces in $\R^4$. Each minimal
surface in $S^3$ generates a minimal foliated semi-symmetric
hypersurface in $\R^4$ according to a special construction given
in \cite{GanMil-CR}. We illustrate how this construction can be
applied  to two examples of minimal surfaces in $S^3$  for
obtaining the first and the second type helicoids, which are
special minimal foliated semi-symmetric hypersurfaces. We also
apply the construction to a  class of generalized tori of second
type, which are new surfaces in the literature, and thus we obtain
new minimal foliated semi-symmetric hypersurfaces in $\R^4$.

\section{Generalized tori in $S^3$}\label{S:Tori}

Let $M^2: l = l(u,v), \,\, (u,v) \in \mathcal{D}$ be a surface,
parameterized by isothermal parameters, and lying on the unit
sphere $S^3$ in $\R^4$,  i.e. the vector-valued function $l(u,v)$
satisfies the equalities:
$$\langle l_u, l_u \rangle = \langle l_v, l_v \rangle = E(u,v); \quad \langle l_u, l_v \rangle = 0; \quad \langle l, l \rangle = 1. \leqno{(2.1)}$$
Since $l, l_u, l_v$ are mutually orthogonal, there exists a unique
(up to a sign) unit vector field $n(u,v)$, such that $\{l, l_u,
l_v, n\}$ form an orthogonal basis in $\R^4$. Differentiating the
equalities (2.1),  we get the following derivative formulas:
$$\begin{array}{l}
\vspace{2mm}
  l_{uu} = \displaystyle{\frac{E_u}{2 E}\, l_u - \frac{E_v}{2 E}\, l_v - E\, l \,\,\, + a_{11} \,n};\\
\vspace{2mm}
l_{uv} = \displaystyle{\frac{E_v}{2 E}\, l_u + \frac{E_u}{2 E}\, l_v \,\,\,\,\,\,\,\quad \quad + a_{12} \,n};\\
\vspace{2mm}
  l_{vv} = \displaystyle{- \frac{E_u}{2 E}\, l_u + \frac{E_v}{2 E}\, l_v - E\, l + a_{22} \,n},
\end{array}$$
where $a_{ij}(u,v), \,\, i,j = 1,2$ are functions defined in
$\mathcal{D}$. Hence,
$$l_{uu} + l_{vv} + 2E\,l = (a_{11}+ a_{22})\, n. \leqno{(2.2)}$$

$M^2$ is a minimal surfaces in $S^3$ if and only if $a_{11}+
a_{22} = 0$. Equality (2.2) implies that  $M^2$ is minimal if and
only if $l_{uu} + l_{vv} + 2E\,l = 0$.

Consequently, the problem of finding the minimal surfaces in $S^3$
is equivalent to the  solvability of the following system
$$\begin{array}{l}
\vspace{2mm}
\Delta l(u,v) + 2 E(u,v)\,l(u,v) = 0; \\
\vspace{2mm} \langle l_u, l_u \rangle = \langle l_v, l_v \rangle =
E(u,v); \quad \langle l_u, l_v \rangle = 0; \quad \langle l, l
\rangle = 1
\end{array} \leqno(2.3)$$
with an appropriate $\mathcal{C}^{\infty}$ smooth scalar function
$E(u,v) >0$ in a small neighbourhood $\mathcal{D}$ of the origin.

Let us note that according to the theorem of H\'{e}lein (see
\cite{G-M-S}, p. 346) the solutions of system (1.1) (and hence of
(2.3)) are $\mathcal{C}^{\infty}$ smooth, because $\mathcal{D}$ is
a two-dimensional domain.

\vskip 1mm Now let $M^2: l = l(u,v)$ be a minimal surface in
$S^3$. Then the derivative formulas of $M^2$ are as follows:
$$\begin{array}{l}
\vspace{2mm}
  l_{uu} = \displaystyle{\frac{E_u}{2 E}\, l_u - \frac{E_v}{2 E}\, l_v - E\, l \,\,\, + a \,n};\\
\vspace{2mm}
l_{uv} = \displaystyle{\frac{E_v}{2 E}\, l_u + \frac{E_u}{2 E}\, l_v \,\,\,\,\,\,\,\quad \quad + b \,n};\\
\vspace{2mm}
  l_{vv} = \displaystyle{- \frac{E_u}{2 E}\, l_u + \frac{E_v}{2 E}\, l_v - E\, l - a
  \,n},
\end{array}\leqno{(2.4)}$$
where $a(u,v)$ and $b(u,v)$ are functions in $\mathcal{D}$.  The
derivatives $n_u$ and $n_v$ of $n(u,v)$ satisfy
$$n_u = \displaystyle{- \frac{a}{E}\, l_u - \frac{b}{E}\,
l_v};\qquad n_v = \displaystyle{- \frac{b}{E}\, l_u +
\frac{a}{E}\, l_v}. \leqno{(2.5)}$$

Using the Gauss and Codazzi equations (or equivalently the
identities of the mixed third derivatives of $l(u,v)$ and mixed
second derivatives of $n(u,v)$) from (2.4) and (2.5) we get that
the functions $a(u,v)$ and $b(u,v)$ are harmonic functions,
satisfying the Cauchy - Riemann conditions:
$$b_u(u,v) = a_v(u,v); \qquad b_v(u,v) = - a_u(u,v).$$
The Gauss and Codazzi equations for $M^2$ also imply the identity
$$a^2 + b^2 = \displaystyle{\frac{1}{2} \Delta E - \frac{E_u^2 + E_v^2}{2 E} + E^2}. \leqno(2.6)$$

As it is well known the Gauss curvature $K$ of $M^2$ is given by
$$K = \displaystyle{\frac{E_u^2 + E_v^2}{2 E^3} - \frac{\Delta E}{2 E^2} =
- \frac{1}{2 E} \Delta \ln E}. \leqno(2.7)$$ Hence, equalities
(2.6) and (2.7) imply that the Gauss curvature $K$ is expressed by
the functions $a$ and $b$ as follows:
$$K = \ds{1 - \frac{a^2 + b^2}{E^2}}.  \leqno(2.8)$$

\vskip 2mm The simplest case, in which problem (2.3) can be solved
completely, is the case $K = const$. In \cite{Lawson-2} it is
proved that: \textit{If $M^2$ is a minimal surface in $S^3$ of
constant Gauss curvature $K$, then either $K = 1$ and $M^2$ is
totally geodesic, or $K = 0$ and $M^2$ is an open piece of the
Clifford torus}.

From (2.8) it follows that the case $M^2$ is totally geodesic in
$S^3$ (i.e. $M^2$ is a sphere with radius 1) corresponds to $a = b
= 0$. Further on we shall consider only the case $(a,b) \neq
(0,0)$.

\vskip 2mm Let us recall that the Clifford torus is a surface  in
$\R^4$, parameterized as follows:
$$\mathcal{M}:  l(u,v) = (\cos u \,\cos v; \cos u\, \sin v; \sin u \,\cos v; \sin u\,\sin v).$$
A direct computation shows that $l(u,v)$ satisfies the equality
$l_{uu} + l_{vv} + 2 l = 0$, and $\langle l_u, l_u \rangle =
\langle l_v, l_v \rangle = 1; \,\, \langle l_u, l_v \rangle = 0;
\,\, \langle l, l \rangle = 1$. The parametric lines $u = const$
and $v = const$ of $\mathcal{M}$ are circles.

Another example of minimal surfaces in $S^3$, which is a
generalization of the Clifford torus,  is given in \cite{Lawson}.
H. B. Lawson proved that every "ruled" minimal surface in $S^3$ is
an open submanifold of one of the surfaces $\mathcal{M}_{\alpha}$
given by
$$\mathcal{M}_{\alpha}:  l(x,y) = (\cos x \,\cos \alpha y; \cos x\, \sin \alpha y; \sin x \,\cos y; \sin x\,\sin y) \leqno(2.9)$$
for some constant $\alpha >0$. Here "ruled" means one-parameter
family of great circles in $S^3$. The surface $\mathcal{M}_1$ is
the Clifford torus, and it is the only surface
$\mathcal{M}_{\alpha}$ with constant Gauss curvature (see
\cite{Lawson-2}). We shall call the surfaces
$\mathcal{M}_{\alpha}$  ($\alpha \neq 1$) \textit{Lawson tori}.

The tangent space of $\mathcal{M}_{\alpha}$ is spanned by the
vector fields
$$\begin{array}{l}
\vspace{2mm}
l_x (x,y) = ( - \sin x\, \cos \alpha y; - \sin x\,\sin \alpha y;  \cos x\,\cos y; \cos x\,\sin y),\\
\vspace{2mm}
l_y (x,y) = ( - \alpha \cos x \,\sin \alpha y; \alpha \cos x\,\cos \alpha y; - \sin x\,\sin y; \sin x\,\cos y);\\
\end{array} \leqno{(2.10)}$$
and the coefficients $E$, $F$, $G$ of the first fundamental form
of $\mathcal{M}_{\alpha}$ are given by $E = 1; \, F = 0; \, G =
G(x) = \alpha^2 \cos^2 x + \sin^2 x$.

Using (2.9) and (2.10) we find the  unit normal vector field
$n(u,v)$ of $\mathcal{M}_{\alpha}$ orthogonal to $\{l, l_u, l_v\}$
:
$$n(x,y) = \ds{\frac{1}{\sqrt{\alpha^2 \cos^2 x + \sin^2 x}}} (\sin x\,\sin \alpha y; - \sin x\,\cos \alpha y; - \alpha \cos x\,\sin y; \alpha \cos x\,\cos y).$$
A direct computation shows that the Gauss curvature of
$\mathcal{M}_{\alpha}$ is given by
$$K = 1 - \ds{\frac{\alpha^2}{(\alpha^2 \cos^2 x + \sin^2 x)^2}}$$
and obviously  $K \neq const$ when  $\alpha \neq 1$.

Let us consider the Lawson torus $\mathcal{M}_{\alpha}$ for
$\alpha \neq 1$. In such case the parametric lines $y = y_0 =
const$ of $\mathcal{M}_{\alpha}$ are circles, while the parametric
lines $x = x_0 = const$ are curves in $\R^4$ with constant Frenet
curvatures.

The parametrization (2.9) of $\mathcal{M}_{\alpha}$ is not an
isothermal one. We shall find isothermal parameters of
$\mathcal{M}_{\alpha}$, which are also principal parameters of the
surface. Let us consider the following change of the parameters:
$$\overline{u} = \ds{\int_0^x \frac{1}{\sqrt{G(\tau)}}\, d\tau}; \qquad \overline{v} = v.$$
Then we obtain
$$\overline{E} = \overline{G} = G(x(\overline{u})); \qquad \overline{F} = 0,$$
i.e. the parameters $(\overline{u},\overline{v})$ are isothermal.
A direct computation shows that the vector-function
$l(\overline{u},\overline{v})$ satisfies the system:
$$l_{\overline{u} \,\overline{u}} = \displaystyle{\frac{\overline{E}_{\overline{u}}}{2 \overline{E}}\, l_{\overline{u}} - \overline{E}\, l};
\qquad
l_{\overline{u} \,\overline{v}} = \displaystyle{\frac{\overline{E}_{\overline{u}}}{2 \overline{E}}\, l_{\overline{v}} + \alpha \,n};\\
\qquad
  l_{\overline{v} \,\overline{v}} = \displaystyle{- \frac{\overline{E}_{\overline{u}}}{2 \overline{E}}\, l_{\overline{u}} - \overline{E}\, l}.$$
Hence, for the Lawson torus $\mathcal{M}_{\alpha}$ the functions
$a$ and $b$ in formulas (2.4) are $a = 0,\,\, b = \alpha = const$,
and $\overline{E} = \overline{E}(\overline{u})$. The circles on
$\mathcal{M}_{\alpha}$ are the parametric $\overline{u}$-lines.

If we change the isothermal parameters
$(\overline{u},\overline{v})$ with isothermal parameters $(u,v)$
in the following way
$$u = \ds{\frac{\sqrt{2}}{2}(\overline{u} + \overline{v})}; \qquad v = \ds{\frac{\sqrt{2}}{2}(\overline{u} - \overline{v})},$$
then $\widetilde{E}(u,v) = \langle l_u, l_u \rangle = \langle l_v,
l_v \rangle  = \widetilde{G}(u,v)$, and $l(u,v)$ satisfies (2.4) with $a = \alpha$, $b =0$, and $\widetilde{E}$ instead of $E$.

With respect to $(u,v)$ the surface $\mathcal{M}_{\alpha}$ is
parameterized by principal lines, i.e. the shape operator,
corresponding to the normal vector field $n(u,v)$, is in diagonal
form. The circles on $\mathcal{M}_{\alpha}$ are bisectrices of the
principal lines.

\vskip 1mm Now we shall find all minimal surfaces in $S^3$,  which
are one-parameter families of circles, crossing  the principal
lines under a constant angle. We call these surfaces
\textit{generalized tori} in $S^3$. They generalize the Lawson
tori.

\begin{prop} \label{P:principal parameters}
Let $M^2$ be a minimal surface in $S^3$ with non-constant Gauss
curvature. Then $M^2$ can locally be parameterized by principal
lines, and the new parameters are isothermal ones.
\end{prop}

\noindent \textit{Proof.} Let $M^2: l = l(u,v), \,\, (u,v) \in
\mathcal{D}$ be a minimal surface in $S^3$, parameterized by
isothermal parameters. Then the derivative formulas (2.4) of $M^2$
hold good, and the functions $a(u,v)$ and $b(u,v)$ are harmonic
functions, satisfying the Cauchy - Riemann conditions. In case of
$b(u,v) \equiv 0$ the parameters $(u,v)$ are principal ones. Let
$b(u_0,v_0) \neq 0$, $(u_0,v_0) \in \mathcal{D}$. Then there
exists $\mathcal{D}_0 \subset \mathcal{D}$ such that $b(u,v) \neq
0$ for all $(u,v) \in \mathcal{D}_0 $. We shall prove that there
exist isothermal parameters $(x,y)$ such that $\overline{b}(x,y) =
\langle l_{x\,y}, n \rangle = 0$. If
$$x = x(u,v); \qquad
y = y (u,v)$$ is  a holomorphic change of the parameters ($x(u,v)$
and $y(u,v)$ satisfy the Cauchy - Riemann conditions), then
$$\overline{b}(x,y) =  2 a u_{x} u_{y} - b (u_{x}^2 - u_{y}^2).$$
Hence,
$$\overline{b}  = 0 \qquad \Longleftrightarrow \qquad
\ds{b \left(\frac{u_{x}}{u_{y}}\right)^2 - 2a
\left(\frac{u_{x}}{u_{y}}\right) - b = 0}.$$

From the inverse change of the parameters, using the Cauchy -
Riemann conditions, we have $\ds{x_u = \frac{u_{x}}{u_{x}^2 +
u_{y}^2}}$,
 $\ds{x_v = - \frac{u_{y}}{u_{x}^2 + u_{y}^2}}$, and hence we obtain that
 $\overline{b}  = 0$ if and only if
$$\ds{b \left(\frac{x_{u}}{x_{v}}\right)^2 + 2a \left(\frac{x_{u}}{x_{v}}\right) - b = 0},$$
i.e. $\ds{\frac{x_u}{x_v} = \frac{-a \pm \sqrt{a^2 + b^2}}{b}}$.
We denote $\ds{\beta(u,v) = \frac{- a + \sqrt{a^2 + b^2}}{b}, \;
\gamma(u,v) = \frac{- a - \sqrt{a^2 + b^2}}{b}}$.

Now, let us consider the equations
$$\ds{\frac{dv}{du} = \beta(u,v); \qquad \frac{dv}{du} = \gamma(u,v)}. \leqno{(2.11)}$$
For each point $(u_0, v_0) \in \mathcal{D}_0$ there exists
$\mathcal{D}_1 \subset \mathcal{D}_0$ and functions $\Phi(u,v)
\neq 0$, $\Psi(u,v) \neq 0$ in $\mathcal{D}_1$, such that the
integral curves of the first equation in (2.11) are given by
$\Phi(u,v) = const$, while the integral curves of the second
equation in (2.11) are $\Psi(u,v)  = const$. Hence, $$\Phi_u =
-\beta \Phi_v, \qquad \Psi_u = -\gamma \Psi_v. \leqno{(2.12)}$$ We
consider the following smooth change of the parameters:
$$x = \Phi(u,v), \qquad y = \Psi(u,v); \qquad (u,v) \in \mathcal{D}_1. \leqno{(2.13)}$$
When $(u,v)$ runs in $\mathcal{D}_1$ the parameters $(x,y)$
describe a domain $\overline{\mathcal{D}} \subset \R^2$. Now $x_u
= \Phi_u; \; x_v = \Phi_v; \; y_u = \Psi_u; \; y_v = \Psi_v$.
Using that $\beta \gamma = -1$ and equalities (2.12), we get
$\langle l_x, \, l_y \rangle = 0$, i.e. the parametrization
$(x,y)$ is orthogonal one. We shall prove that this
parametrization is isothermal. With respect to the new parameters
the coefficients of the first fundamental form are $E = \langle
l_x, l_x \rangle$, $F = \langle l_x, l_y \rangle = 0$, and $G =
\langle l_y, l_y \rangle$. Then for the surface $M^2$ the
following derivative formulas hold:
$$\begin{array}{l}
\vspace{2mm}
  l_{xx} = \displaystyle{\frac{E_x}{2 E}\, l_x - \frac{E_y}{2 E}\, l_y - E\, l \,\,\, + a \,n};\\
\vspace{2mm}
l_{xy} = \displaystyle{\frac{E_y}{2 E}\, l_x + \frac{G_x}{2 G}\, l_y};\\
\vspace{2mm}
  l_{yy} = \displaystyle{- \frac{G_x}{2 G}\, l_x + \frac{G_y}{2 G}\, l_y - G\, l - a \,n},
\end{array} \leqno{(2.14)}$$

$$n_x = \displaystyle{- \frac{a}{E}\, l_x };\qquad
n_y = \displaystyle{\frac{a}{G}\, l_y}.
\leqno{(2.15)}$$
Taking into account the second fundamental form of $M^2$ as a
surface in $\R^4$, from (2.14) we calculate that the Gauss
curvature $K$ is given by
$$K = \ds{1 - \frac{a^2}{E G}}. \leqno{(2.16)}$$
 Using that $n_{xy} = n_{yx}$, $l_{xxy} = l_{xyx}$, $l_{xyy}
= l_{yyx}$, from (2.14) and (2.15) we obtain
$$a_x = 0; \qquad a_y = 0; \qquad (E - G)E_y =0; \qquad (E - G) G_x = 0;\leqno{(2.17)}$$
$$\begin{array}{l}
\vspace{2mm}
\ds{\frac{E_{yy}}{2E} + \frac{G_{xx}}{2G} - \frac{3E_y^2}{4E^2} - \frac{G_x^2}{4G^2} - \frac{E_x G_x - E_y G_y}{4EG} + E - \frac{a^2}{G} = 0};\\
\vspace{2mm} \ds{\frac{E_{yy}}{2E} + \frac{G_{xx}}{2G} -
\frac{3G_x^2}{4G^2} - \frac{E_y^2}{4E^2} + \frac{E_x G_x - E_y
G_y}{4EG} + G - \frac{a^2}{E} = 0}.
\end{array}  \leqno{(2.18)}$$

If we assume that $E(x_0,y_0) - G(x_0,y_0) \neq 0$ at some point
$(x_0,y_0) \in \overline{\mathcal{D}}$, and hence $E - G \neq 0$
in a neighbourhood $\overline{\mathcal{D}}_0 \subset
\overline{\mathcal{D}}$ of $(x_0,y_0)$, then from (2.17) we get
$E_y = G_x = 0$ in $\overline{\mathcal{D}}_0$. Now equalities
(2.16), (2.17), and (2.18) imply that $K = 0$ in
$\overline{\mathcal{D}}_0$, which contradicts the assumption in
the theorem. Hence, $E - G \equiv 0$ in $\overline{\mathcal{D}}$.
Consequently, the parameters $(x,y)$, defined by change (2.13),
are principal ones. \qed

\vskip 3mm Now let $M^2: l = l(u,v)$ be a minimal surface in
$S^3$, parameterized locally by isothermal principal parameters,
i.e. $b = 0$. Using that $b_u = a_v, \,\, b_v = - a_u$, we get $a
= const$. Without loss of generality we assume that $a = 1$ (if $a
\neq 1$, we multiply the parameters by $\sqrt{| a|}$). Hence, the
derivative formulas (2.4) and (2.5) hold with $a = 1$, $b = 0$.

\begin{thm} \label{T:generalized torus}
Let $M^2$ be a minimal surface in $S^3$ with non-constant Gauss
curvature. If on $M^2$ there exists a family of circles, crossing
the principal lines under a constant angle $\theta $, then the
circles are either principal lines ($\theta  = 0$ or $\theta  =
\ds{\frac{\pi}{2}}$) or  bisectrices of the principal lines
($\theta  = \ds{\frac{\pi}{4}}$ or $\theta  =
\ds{\frac{3\pi}{4}}$).
\end{thm}

\noindent \textit{Proof.} Let $M^2: l = l(u,v)$ be parameterized
locally by principal parameters. Suppose that on $M^2$ there
exists a family of circles, crossing  the principal lines under a
constant angle $\theta$. Let
$$\begin{array}{l}
\vspace{2mm}
  x = \cos \theta\, u + \sin \theta \,v;\\
\vspace{2mm}
  y = -\sin \theta \,u + \cos \theta \,v;
\end{array}
\qquad \theta = const, \,\, \theta \in [0; 2\pi). \leqno{(2.19)}$$
Then from (2.19) we get
$$E_x = \cos \theta  \,E_u + \sin \theta  \,E_v; \qquad
E_y = -\sin \theta  \,E_u + \cos \theta \,E_v. \leqno{(2.20)}$$
Using (2.4) with $a=1$, $b=0$,  and (2.20) we calculate
$$\begin{array}{l}
\vspace{2mm}
  l_{xx} = \displaystyle{\frac{E_x}{2 E}\, l_x - \frac{E_y}{2 E}\, l_y - E\, l + \cos 2\theta \,n};\\
\vspace{2mm}
l_{xy} = \displaystyle{\frac{E_y}{2 E}\, l_x + \frac{E_x}{2 E}\, l_y} \;\;\;\;\;\;\;\;\;\; - \sin 2\theta  \,n;\\
\vspace{2mm}
  l_{yy} = \displaystyle{- \frac{E_x}{2 E}\, l_x + \frac{E_y}{2 E}\, l_y -E\, l - \cos 2\theta \,n}.
\end{array}\leqno{(2.21)}$$
Let us denote $a = \cos 2\theta , \,\, b = \sin 2\theta $ ($a$,
$b$ - const).

Now we shall consider an arbitrary $x$-line $c: l(x) = l(x,y_0)$.
The curve $c$ is a circle if and only if its Frenet curvatures are
$\varkappa = const;\,\, \tau = \sigma = 0$. Using (2.21) we
calculate the tangent vector $t_c$ and the principal normal vector
$n_c$ of $c$:
$$t_c = \ds{\frac{l_x}{\sqrt{E}}}; \qquad
\vspace{2mm} n_c = \ds{\frac{1}{\varkappa} \left(-\frac{E_y}{2E^2}
\,l_y - l + \frac{a}{E}\, n \right)},$$
where $\varkappa^2 = \ds{\frac{E_y^2}{4E^3} + 1 +
\frac{a^2}{E^2}}$. The derivatives of $t_c$ and $n_c$ are:
$$\begin{array}{l}
\vspace{2mm}
t_c' = \varkappa \,n_c;\\
\vspace{2mm} n_c' = \ds{- \varkappa \, t_c -
\left[\left(\frac{1}{\varkappa}\right)' \frac{E_y}{2E^2} +
       \frac{1}{\varkappa} \left( \left(\frac{E_y}{2E^2} \right)_x +
       \frac{E_x E_y}{4E^3} + \frac{a \,b}{E^2}\right)\right] \frac{l_y}{\sqrt{E}}}\\
\qquad \quad - \ds{\left(\frac{1}{\varkappa}\right)' \,
\frac{l}{\sqrt{E}} + \left[\left(\frac{1}{\varkappa}\right)'
\frac{a}{E} +
       \frac{1}{\varkappa} \left( - b \frac{E_y}{2E^2} +
       \left(\frac{a}{E}\right)_{x}\right)\right] \, \frac{n}{\sqrt{E}} }.
\end{array} \leqno{(2.22)}$$
From (2.22) it follows that $c$ is a circle if and only if
$$\begin{array}{l}
\vspace{2mm}
\varkappa = const;\\
\vspace{2mm} \ds{\left(\frac{E_y}{2E^2}
\right)_x + \frac{E_x E_y}{4E^3} + \frac{a \,b}{E^2} = 0};\\
\vspace{2mm}
 \ds{ \left(\frac{a}{E}\right)_{x} - b \frac{E_y}{2E^2} = 0},
\end{array}$$
which is equivalent to
$$\begin{array}{l}
\vspace{2mm}
2a E_x + b E_y = 0;\\
\vspace{2mm} 2E\, E_{xy} - 3 E_x E_y + 4 a\,b\,E = 0.
\end{array} \leqno{(2.23)}$$

Analogously, the $y$-lines are circles  if and only if
$$\begin{array}{l}
\vspace{2mm}
b E_x - 2 a E_y = 0;\\
\vspace{2mm} 2E\, E_{xy} - 3 E_x E_y  - 4 a\,b\,E  = 0.
\end{array} \leqno{(2.24)}$$
From (2.20) we calculate
$$E_{xy} = - \sin \theta \cos \theta E_{uu} + \cos 2\theta E_{uv} + \sin \theta \cos \theta E_{vv}. \leqno{(2.25)}$$

From the first equality of (2.23), using (2.20), we get $\cos^3
\theta \,E_u - \sin^3 \theta \,E_v = 0$. All solutions of this
equation are given by
$$E = \varphi(\sin^3 \theta \, u + \cos^3 \theta \,v)   \leqno{(2.26)}$$
for arbitrary smooth function $\varphi$. Hence, $E_u = \sin^3
\theta \, \varphi'; \,\, E_v = \cos^3 \theta \, \varphi';$ $E_{uu}
= \sin^6 \theta \, \varphi'';$ $E_{vv} = \cos^6 \theta \,
\varphi'';$ $E_{uv} = \sin^3 \theta\, \cos^3 \theta \, \varphi''$.
From the second equality of (2.23), using (2.25) and (2.26) we
obtain $\ds{\sin 2\theta \, \cos 2 \theta (\varphi \varphi'' -
\frac{3}{2}\varphi'^2 - 4 \varphi) = 0}$. Consequently, the
$x$-lines are circles if and only if
$$\begin{array}{l}
\vspace{2mm}
E = \varphi(\sin^3 \theta \, u + \cos^3 \theta \,v);\\
\vspace{2mm} \sin 2\theta \, \cos 2 \theta (\varphi \varphi'' -
\ds{\frac{3}{2}}\varphi'^2 - 4 \varphi) = 0.
\end{array} \leqno{(2.27)}$$
Analogously, using (2.24)  we obtain that the $y$-lines are
circles if and only if
$$\begin{array}{l}
\vspace{2mm}
E = \varphi(\cos^3 \theta \, u - \sin^3 \theta \,v);\\
\vspace{2mm} \sin 2\theta \, \cos 2 \theta (\varphi \varphi'' -
\ds{\frac{3}{2}}\varphi'^2 - 4 \varphi) = 0.
\end{array} \leqno{(2.28)}$$

Thus the condition the $x$-lines (or $y$-lines) to be circles
leads to the following cases:

\vskip 1mm I. $\sin 2\theta = 0$, i.e. $\theta = 0$ or $\theta =
\ds{\frac{\pi}{2}}$.

\noindent This case corresponds to $x = u; \,\, y = v$ or $x = v;
\,\, y = - u$. From (2.27) and (2.28) we obtain that the $u$-lines
are circles if and only if $E = \varphi(v)$, and the $v$-lines are
circles if and only if $E = \varphi(u)$. In this case one of the
families of principal lines is a family of circles.

\vskip 1mm II. $\cos 2\theta = 0$, i.e. $\theta =
\ds{\frac{\pi}{4}}$ or $\theta = \ds{\frac{3\pi}{4}}$.

\noindent This case corresponds to  $x = \frac{\sqrt{2}}{2}(u +
v); \,\, y = \frac{\sqrt{2}}{2}(- u + v)$ or $x =
\frac{\sqrt{2}}{2}( - u + v); \,\, y = - \frac{ \sqrt{2}}{2}( u +
v)$. From (2.27) and (2.28) we obtain that the $x$-lines are
circles if and only if $E = \varphi(x)$, and the $y$-lines are
circles if and only if $E = \varphi(y)$. In this case one of the
families of bisectrices of the principal lines is a family of
circles.

\vskip 1mm III. $\varphi \varphi'' - \ds{\frac{3}{2}}\varphi'^2 -
4 \varphi = 0$, and $\sin 2\theta\,\cos 2\theta  \neq 0$.

\noindent We shall prove that this case is not possible. Since $E
> 0$ then $E = \varphi(\tau) = e^{z(\tau)}$ for some function $z =
z(\tau)$, $\tau = \sin^3 \theta \, u + \cos^3 \theta \,v$ (or
$\tau = \cos^3 \theta \, u - \sin^3 \theta \,v$). Moreover,  $E
\neq const$, i.e. $z'(\tau) \neq 0$. The equality $\varphi
\varphi'' - \ds{\frac{3}{2}}\varphi'^2 - 4 \varphi = 0$ implies
$$z'' - \frac{1}{2} z'^2 - 4 e^{-z} = 0. \leqno(2.29)$$

On the other hand, using identity (2.6), we obtain
$$(\cos^6 \theta + \sin^6 \theta) z'' + 4 \sinh z = 0. \leqno(2.30)$$
Let us denote $\lambda = \cos^6 \theta + \sin^6 \theta = const$.
Multiplying (2.30) by $z'$ and integrating we get
$$\frac{\lambda}{2} z'^2 (\tau) + 4 \cosh z(\tau) = \frac{\lambda}{2} z'^2 (0) + 4 \cosh z(0).$$
Equalities (2.29) and (2.30) imply
$$4 \sinh z + \frac{\lambda}{2} z'^2 + 4 \lambda e^{-z} = 0.$$
Using the last two equalities we obtain
$$4(1 - \lambda) e^{-z} = \ds{\frac{\lambda}{2} z'^2 (0) + 4 \cosh z(0) =
const.}$$
 Since $1 - \lambda = 3 \sin^2 \theta \, \cos^2 \theta
\neq 0$, we get $e^{-z(\tau)} = const$, i.e. $z(\tau) = const$,
which contradicts the condition $z'(\tau) \neq 0$, i.e. $E \neq
const$. \qed

\vskip 3mm From Theorem \ref{T:generalized torus} it follows that
there are only two types of generalized tori in $S^3$: the first
one is characterized by the condition that one of the families of
bisectrices of the  principal lines is a family  of circles (we
shall call such surfaces \textit{generalized tori of first type}),
and the second one is characterized by the condition that one of
the families of principal lines is a family of circles (we shall
call these surfaces \textit{generalized tori of second type}).

\begin{thm} \label{T:First type}
Let $M^2$ be a generalized torus of first type with non-constant
Gauss curvature. Then $M^2$ is a Lawson torus
$\mathcal{M}_{\alpha}$ for some $\alpha >0$, $\alpha \neq 1$.
\end{thm}

\noindent \textit{Proof}. Let $M^2$ be a generalized torus of
first type with non-constant Gauss curvature. In this case the
derivative formulas (2.4) hold, with $a = 0;\,\, b = 1$, and $E =
E(u)$ (or $E = E(v)$). We shall consider only the first case, i.e.
$E = E(u)$, $E_u \neq 0$. The second one can be investigated
analogously. The derivative formulas in this case look like:
$$\begin{array}{ll}
\vspace{2mm}
  \ds{l_{uu} =  \frac{E_u}{2E}\, l_u - E \,l}; \qquad & \ds{ n_u = - \frac{1}{E} \,l_v};\\
\vspace{2mm}
\ds{l_{uv} =  \frac{E_u}{2E}\, l_v + n}; \qquad & \ds{n_v = - \frac{1}{E} \,l_u};\\
\vspace{2mm}
   \ds{l_{vv} = -\frac{E_u}{2E}\, l_u  -E \,l};
\end{array} \leqno{(2.31)}$$

We shall prove that the parametric $u$-lines are great circles.
Let $c: c(u) = l(u,v_0)$, $v_0 = const$ be an arbitrary $u$-line.
From (2.31) it follows that $\dot{c} = l_u$, and the tangent
vector $t_c$ of $c$ is $t_c = c' = \ds{\frac{\dot{c}}{\dot{s}} =
\frac{l_u}{\sqrt{E}}}$. We calculate $t_c' =
\ds{\frac{\dot{t}_c}{\dot{s}} = \frac{1}{\sqrt{E}}
\left(\frac{l_{uu}}{\sqrt{E}} - \frac{E_u}{2E \sqrt{E}}\,l_u
\right) = -l}$. Hence, the curvature $\varkappa$ of $c$ is
$\varkappa = 1$, i.e. $c$ is a great circle. Consequently, $M^2$
is a one-parameter family of great circles. According to \cite[
Prop. 7.2]{Lawson} $M^2$ is on open submanifold of
$\mathcal{M}_{\alpha}$ for some $\alpha >0$, $\alpha \neq 1$. \qed

\vskip 3mm Now we shall consider a generalized torus of second
type with non-constant Gauss curvature. In this case the
derivative formulas (2.4) hold, with $a = 1;\,\, b = 0$, and
$E = E(u)$ (or $E = E(v)$). We consider the case $E = E(u)$, $E_u
\neq 0$. In this case the parametric $v$-lines are circles. We
shall prove that the different $v$-lines are circles with
different radius.

Now the derivative formulas are as follows:
$$\begin{array}{ll}
\vspace{2mm}
  \ds{l_{uu} =  \frac{E_u}{2E}\, l_u - E \,l + n}; \qquad & \ds{n_u = - \frac{1}{E} \,l_u};\\
\vspace{2mm}
\ds{l_{uv} =  \frac{E_u}{2E}\, l_v }; \qquad & \ds{n_v =  \frac{1}{E} \,l_v};\\
\vspace{2mm}
   \ds{l_{vv} = -\frac{E_u}{2E}\, l_u  -E \,l - n};
\end{array} \leqno{(2.32)}$$
Let $c: c(v) = l(u_0,v)$, $u_0 = const$ be an arbitrary $v$-line.
As in the proof of Theorem \ref{T:First type}, from (2.32) we
calculate the curvature $\varkappa(u_0)$ of $c$: $\varkappa(u_0) =
\ds{\sqrt{1+ \frac{E_u^2(u_0)}{4 E^3(u_0)} +
\frac{1}{E^2(u_0)}}}$. For different values of the constant $u_0$
the curvatures $\varkappa(u_0)$ are different. We note that
$\varkappa(u_0) >1$, i.e. the circles are not great ones.
Therefore, the generalized tori of second type are different from
the Lawson tori.

\vskip 1mm Since $E(u) > 0$, we  write $E(u)$ in the form $E(u) =
e^{z(u)}$, $z(u) \neq const$. Then the system (2.32) is rewritten
in the following form
$$\begin{array}{ll}
\vspace{2mm}
  \ds{l_{uu} - \frac{z'(u)}{2}\, l_u  + e^{z(u)}\,l - n = 0}; \qquad & \ds{ n_u + e^{- z(u)}\, l_u  = 0};\\
\vspace{2mm}
\ds{l_{uv} - \frac{z'(u)}{2}\, l_v = 0}; \qquad & \ds{n_v - e^{- z(u)}\,l_v = 0};\\
\vspace{2mm}
   \ds{l_{vv} + \frac{z'(u)}{2}\, l_u  + e^{z(u)}\,l + n = 0}.
\end{array} \leqno{(2.33)}$$
We look for classical solutions of the above system for the
vector-valued functions $l(u,v)$, $l_u(u,v)$, $l_v(u,v)$ and
$n(u,v)$ in a neighbourhood of the origin under the following
initial conditions:
$$\begin{array}{ll}
\vspace{2mm}
l(0,0) = e_1; \qquad \qquad \qquad & \ds{l_v(0,0) =  e^{\frac{s}{2}} e_3 = \sqrt{E(0)}\, e_3};\\
\vspace{2mm}
\ds{l_u(0,0) = e^{\frac{s}{2}} e_2 = \sqrt{E(0)}\, e_2}; \qquad \qquad \qquad & n(0,0) = e_4,\\
\end{array} \leqno{(2.34)}$$
where $\{e_1, e_2, e_3, e_4\}$ is the standard orthonormal basis
in $\R^4$, and $s = const = z(0) = \ln E(0)$.

Since the function $E(u)$ satisfies  identity (2.6) with $a =
1,\,\,b = 0$,  it follows that $z(u)$ is the solution of the
following ordinary differential equation
$$z''(u)+  \displaystyle{4 \sinh z(u)} = 0; \qquad z(0) = s; \quad z'(0) = 2t, \leqno{(2.35)}$$
where $s$ and $t$ are arbitrary constants.

In order to find explicitly $z(u)$ we note that the identity (2.6)
holds for an arbitrary minimal surface in $S^3$, parameterized by
isothermal parameters. So, let us consider again the Lawson torus
$\mathcal{M}_{\alpha}$, defined by (2.9). We change the parameters
$(x,y)$ by new parameters $(u,v)$ in the following way \label{key}
$$u = h(x) = \ds{\sqrt{\alpha} \int_0^x \frac{1}{\sqrt{\alpha^2 \cos^2 \tau + \sin ^2 \tau}}\, d\tau}; \qquad v = \sqrt{\alpha} y,$$
and denote by $h^{-1}$ the inverse function of $h$.  Then we obtain
an isothermal parametrization of $\mathcal{M}_{\alpha}$, and the
function
$$\ds{E(u) = e^{z(u)} = \frac{\alpha^2 \cos^2 h^{-1}(u) + \sin ^2 h^{-1}(u)}{\alpha}}$$
satisfies (2.6) with $a = 0;\,\, b = 1$. Hence the function
$$z(u) = \ds{ \ln \frac{\alpha^2 \cos^2  h^{-1}(u) + \sin ^2  h^{-1}(u)}{\alpha}} \leqno{(2.36)}$$ is a solution of
(2.35) with $t = 0$ and $s = \ln \alpha$, $\alpha > 0, \,\, \alpha
\neq 1$. (It can be calculated directly that the function $z(u)$,
defined by (2.36) satisfies (2.35)  with $z(0) = \ln \alpha$,
$\alpha
> 0, \,\, \alpha \neq 1$ and $z'(0) = 0$.)

We will prove that every solution $\widetilde{z}(u)$ of (2.35)
with arbitrary $t$ and $s$ can be obtained from (2.36) by the
formula $\widetilde{z}(u) = z(u + u_0)$ for a suitable choice of
constants $u_0$ and $\alpha$. Since (2.35) is an autonomous
equation and $z(u)$ is its solution, then it follows that
$\widetilde{z}(u) = z(u + u_0)$ is also a solution of this
equation. Therefore we have to check only the initial conditions.
Let $x_0 = h^{-1}(u_0)$. Simple calculations give us the
equalities
$$\widetilde{z}(0) = z(u_0) = \ds{ \ln \frac{\alpha^2 \cos^2 x_0 + \sin ^2 x_0}{\alpha}} = s;$$
$$\widetilde{z}'(0) = z'(u_0) = \ds{\frac{(1-\alpha^2) \sin 2 x_0}{\sqrt{\alpha} \sqrt{\alpha^2 \cos^2 x_0 + \sin^2 x_0}}} = 2t,$$
which imply
$$\sin 2 x_0 = \ds{\frac{2 \alpha t e^{\frac{s}{2}}}{1 - \alpha^2}}; \qquad \cos 2 x_0 = \ds{\frac{1  + \alpha^2 - 2 \alpha e^s}{1 - \alpha^2}}.$$
Using that $\sin^2 2 x_0 + \cos^2 2 x_0 = 1$, we get the following
equation for $\alpha$:
$$e^s \alpha^2 - (1+ e^{2s} + t^2 e^s) \alpha + e^s = 0,$$
whose positive solutions are
$$\alpha = \ds{\frac{1+ e^{2s} + t^2 e^s \pm \sqrt{(1+ e^{2s} + t^2 e^s)^2 - 4 e^{2s}}}{2 e^s}}.$$
For the above choice of $\alpha$ and $u_0 = h(x_0) =
\ds{\frac{1}{2} \arctan \frac{2 \alpha t e^{\frac{s}{2}}}{1 +
\alpha^2 - 2 \alpha e^s }}$ the initial conditions are satisfied.

It is curious to mention that all solutions of (2.35) are periodic
ones with period $\omega = h(\pi)$.

In order to simplify system (2.33) we change the vector-function
$l(u,v)$ with vector-function $L(u,v)$ determined by
$$l(u,v) = e^{\frac{z(u)}{2}} L(u,v). \leqno{(2.37)}$$
We get the system
$$\begin{array}{ll}
\vspace{2mm}
L_{uu}(u,v) + \ds{ \frac{z'(u)}{2} L_u(u,v) + e^{-z(u)} L(u,v) - e^{-\frac{z(u)}{2}} n = 0};\\
\vspace{2mm}
L_{uv}(u,v) =0;\\
\vspace{2mm}
L_{vv}(u,v) + \ds{ \frac{z'(u)}{2} L_u(u,v) + \left[\left(\frac{z'(u)}{2}\right)^2 + e^{z(u)} \right]L(u,v) + e^{-\frac{z(u)}{2}} n = 0};\\
\vspace{2mm}
n_u + \ds{e^{-\frac{z(u)}{2}}\left(L_u + \frac{z'(u)}{2} L \right) = 0};\\
\vspace{2mm} n_v - e^{-\frac{z(u)}{2}} L_v = 0.
\end{array} \leqno{(2.38)}$$
The initial conditions (2.34) for $l(u,v)$ imply the following
initial conditions for $L(u,v)$:
$$\begin{array}{ll}
\vspace{2mm}
L(0,0) = e^{- \frac{s}{2}}e_1; \qquad \qquad \qquad & \ds{L_v(0,0) =  e_3};\\
\vspace{2mm} \ds{L_u(0,0) = - t e^{-\frac{s}{2}} e_1 + e_2};
\qquad \qquad \qquad & n(0,0) = e_4.
\end{array}$$

From the second equality in (2.38) it follows that $L(u,v) = f(u)
+ g(v)$, where $f(u)$ and $g(v)$ are vector-functions, satisfying
the system
$$\begin{array}{ll}
\vspace{2mm}
f''(u) +  \ds{\frac{z'(u)}{2} f'(u) + e^{-z(u)} (f(u) + g(v)) - e^{-\frac{z(u)}{2}} n = 0};\\
\vspace{2mm}
g''(v) + \ds{ \frac{z'(u)}{2} f'(u) + \left[\left(\frac{z'(u)}{2}\right)^2 + e^{z(u)} \right] (f(u) + g(v)) + e^{-\frac{z(u)}{2}} n = 0};\\
\vspace{2mm}
\ds{n_u + e^{-\frac{z(u)}{2}}\left(f'(u) + \frac{z'(u)}{2} (f(u) + g(v)) \right) = 0};\\
\vspace{2mm}
\ds{n_v - e^{-\frac{z(u)}{2}} g'(v) = 0};\\
\vspace{2mm} \ds{f'(0) = - t e^{-\frac{s}{2}} e_1 + e_2; \qquad
g'(0) = e_3; \qquad f(0) + g(0) = e^{- \frac{s}{2}}e_1}.
\end{array} \leqno{(2.39)}$$
Without loss of generality we assume that $g(0) = 0;\,\, f(0) =
e^{- \frac{s}{2}}e_1$.

Let us fix $u = 0$ in the fourth equality of (2.39). Then after
integration we get
$$n(0,v) = e^{- \frac{s}{2}} g(v) + e_4.$$
Now using the second equality of (2.39) we obtain that $g(v)$
satisfies the initial problem
$$\begin{array}{ll}
\vspace{2mm}
g''(v) + (t^2 + 2 \cosh s) g(v) = - e^{\frac{s}{2}} e_1 - t e_2 - e^{-\frac{s}{2}} e_4;\\
\vspace{2mm} g(0) = 0; \quad g'(0) = e_3.
\end{array}$$
Simple computations give us
$$\begin{array}{ll}
\vspace{2mm} g(v) = & \ds{\frac{1}{t^2 + 2 \cosh s}  (\cos
\sqrt{t^2 + 2 \cosh s} \,v - 1)
(e^{\frac{s}{2}} e_1 + t e_2 + e^{-\frac{s}{2}} e_4)} +\\
\vspace{2mm}
 & \ds{ \frac{1}{\sqrt{t^2 + 2 \cosh s}} \sin \sqrt{t^2 + 2 \cosh s} \,v \, e_3}.
 \end{array} \leqno{(2.40)}$$

Now, multiplying (2.35) with $z'(u)$ and integrating from $0$ to
$u$, we get that $z'(u)$ satisfies the equality
$$(z'(u))^2 + 8 \cosh z(u) = 4 t^2 + 8 \cosh s. \leqno{(2.41)}$$
Using the first and the second equality of (2.39), (2.40) and
(2.41), setting $v = 0$ we obtain that $f(u)$ satisfies the
initial problem:
$$\begin{array}{ll}
\vspace{2mm}
f''(u) +  z'(u) f'(u) + (t^2 + 2 \cosh s)f(u) = e^{\frac{s}{2}}e_1 + t e_2 +  e^{-\frac{s}{2}} e_4;\\
\vspace{2mm} f(0) = e^{-\frac{s}{2}} e_1; \qquad f'(0) = - t
e^{-\frac{s}{2}} e_1 + e_2.
\end{array}$$
Therefore the solution $f(u)$ of the above system can be written
in the form
$$f(u) = \ds{p(u) + \frac{1}{t^2 + 2 \cosh s} (e^{\frac{s}{2}}e_1 + t e_2 +  e^{-\frac{s}{2}} e_4)},$$
where $p(u)$ is the unique solution of the following linear
homogenous system:
$$\begin{array}{ll}
\vspace{2mm}
p''(u) +  z'(u) p'(u) + (t^2 + 2 \cosh s)p(u) = 0;\\
\vspace{2mm} p(0) = \ds{\frac{1}{t^2 + 2 \cosh s}\left[
e^{-\frac{s}{2}}(t^2 + e^{-s}) e_1 - t e_2 -  e^{-\frac{s}{2}}
e_4\right]}; \qquad p'(0) = - t e^{-\frac{s}{2}} e_1 + e_2.
\end{array}$$

Now, using (2.37) we obtain that the solution $l(u,v)$ of problem
(2.33), (2.34) is given by
$$\begin{array}{ll}
\vspace{2mm} l(u,v) = & \ds{e^{\frac{z(u)}{2}} p(u) +
\frac{e^{\frac{z(u)}{2}}}{t^2 + 2 \cosh s} \cos \sqrt{t^2 + 2
\cosh s} \,v
\left(e^{\frac{s}{2}} e_1 + t e_2 + e^{-\frac{s}{2}} e_4\right)} +\\
\vspace{2mm}
 & \ds{ \frac{e^{\frac{z(u)}{2}}}{\sqrt{t^2 + 2 \cosh s}} \,\sin \sqrt{t^2 + 2 \cosh s} \,v \,\, e_3}.
 \end{array}$$

If we denote $\beta = \sqrt{t^2 + 2 \cosh s}$  ($\beta = const$),
then $l(u,v)$ is rewritten in the form:
$$l(u,v) = \ds{e^{\frac{z(u)}{2}} p(u) +
  \frac{e^{\frac{z(u)}{2}}}{\beta^2}\left[ \cos \beta\, v \left(e^{\frac{s}{2}} e_1 + t e_2 + e^{-\frac{s}{2}} e_4\right) +
 \beta \sin \beta\, v \,\,e_3\right]}.$$
We remind that the function $z(u)$ is given explicitly  by (2.36).

Thus we prove the following result:

\vskip 2mm
\begin{thm} \label{T:Second type}
Let $M^2: l = l(u,v)$ be a generalized torus of second type with
non-constant Gauss curvature. Then
$$l(u,v) = \ds{e^{\frac{z(u)}{2}} p(u) +
  \frac{e^{\frac{z(u)}{2}}}{\beta^2}\left[ \cos \beta\, v \left(e^{\frac{s}{2}} e_1 + t e_2 + e^{-\frac{s}{2}} e_4\right) +
 \beta \sin \beta\, v \,\,e_3\right]},$$
where $\beta = \sqrt{t^2 + 2 \cosh s}$ for arbitrary constants $s$
and $t$, the scalar function $z(u)$ is the solution of the
equation
$$z'' + 4 \sinh z(u) = 0; \qquad z(0) = s; \quad z'(0) = 2 t,$$
and is given explicitly by $(2.36)$, $p(u)$ is a vector function,
which is a solution of the system
$$\begin{array}{ll}
\vspace{2mm}
p''(u) +  z'(u) p'(u) + \beta^2 p(u) = 0;\\
\vspace{2mm} p(0) = \ds{\frac{1}{\beta^2}\left[
e^{-\frac{s}{2}}(t^2 + e^{-s}) e_1 - t e_2 -  e^{-\frac{s}{2}}
e_4\right]}; \qquad p'(0) = - t e^{-\frac{s}{2}} e_1 + e_2,
\end{array} \leqno{(2.42)}$$
and $e_1, e_2, e_3, e_4$ is the standard orthonormal basis in
$\R^4$.
\end{thm}

\section{Application to the theory of minimal foliated semi-symmetric hypersurfaces} \label{S:Application}
In this section we shall relate the theory of minimal surfaces in
$S^3$ with the theory of minimal foliated semi-symmetric
hypersurfaces in $\R^4$.

For an $n$-dimensional Riemannian manifold $(M^n,g)$ we denote by
$T_pM^n$ the tangent space to $M^n$ at a point $p \in M^n$ and by
$\X M^n$ - the algebra of all vector fields on $M^n$. The
associated Levi-Civita connection of the metric $g$ is denoted by
$\nabla$, the Riemannian curvature tensor $R$ is defined by
$$R(X,Y) = [\nabla_X, \nabla_Y] -  \nabla_{[X,Y]};
\quad X,Y \in \X M^n.$$

A {\it semi-symmetric space} is a Riemannian manifold $(M^n,g)$,
whose curvature tensor $R$ satisfies the identity
$$ R(X,Y) \cdot R = 0$$
for all vector fields $X, Y \in \X M^n$. (Here $R(X,Y)$ acts as a
derivation on $R$).

According to the classification of Z. Szab\'{o} \cite{Szabo} the
main class of semi-symmetric spaces is the class of all Riemannian
manifolds foliated by Euclidean leaves of codimension two.

The foliated semi-symmetric hypersurfaces in Euclidean space
$\E^{n+1}$ are the hypersurfaces of type number two, i.e.
hypersurfaces whose rank of the second fundamental form is equal
to two everywhere. They are characterized by a second fundamental
form
$$h = \nu_1 \, \eta_1 \otimes \eta_1 + \nu_2 \, \eta_2 \otimes \eta_2, \quad \nu_1 \nu_2 \neq 0,$$
where $\eta_1$ and $\eta_2$ are unit one-forms, $\nu_1$ and
$\nu_2$ are functions on the hypersurface $M^n$. The Euclidean
leaves of the foliation are the integral submanifolds of the
distribution $\Delta_0$, determined by the one-forms $\eta_1$ and
$\eta_2$, i.e. $\Delta_0(p) = \{X \in T_pM^n \,\, | \,\, \eta_1(X)
= 0, \,\, \eta_2(X) = 0\}, \,\, p \in M^n$. A special class of
foliated semi-symmetric hypersurfaces is the class of ruled
hypersurfaces.

A hypersurface $M^n$ of type number two is minimal if $\nu_1 +
\nu_2 = 0$.

\vskip 1mm The foliated semi-symmetric hypersurfaces in $\E^{n+1}$
are characterized in \cite{GanMil-MathB} by the following

\vskip 2mm
\begin{thm}\label{T:envelope}
A hypersurface $M^n$ in Euclidean space $\E^{n+1}$ is locally a
foliated semi-symmetric hypersurface if and only if it is the
envelope of a two-parameter family of hyperplanes in $\E^{n+1}$.
\end{thm}

Using the characterization of a foliated semi-symmetric
hypersurface as the envelope of a two-parameter family of
hyperplanes, each such hypersurface is determined by a pair of a
unit vector-valued function $l(u,v)$ and a scalar function
$r(u,v)$, defined in a domain ${\mathcal D} \subset \R^2$.

Since the vector fields $l_u$ and $l_v$ are linearly independent,
then the vector-valued function $l(u,v)$ determines a
two-dimensional surface $M^2: l = l(u,v)$, $(u,v)\in {\mathcal D}$
in $\E^{n+1}$. Without loss of generality it can be assumed that
the surface $M^2$ is parameterized locally by isothermal
parameters, i.e. $E = G,\, F = 0$. Then the generated foliated
semi-symmetric hypersurface $M^n$ is given in \cite{GanMil-CR} by
$$X(u, v, w^{\alpha}) = r\,l + \displaystyle{\frac{r_u}{E}\,\,l_u} + \displaystyle{\frac{r_v}{E}\,\,l_v} +
w^{\alpha}\,b_{\alpha}, \quad \alpha = 1, \dots, n-2,
\leqno{(3.1)}$$ where $(u,v) \in {\mathcal D}, \,\, w^{\alpha} \in
\R, \,\, \alpha = 1, \dots, n-2$, and $b_1(u,v), \dots,
b_{n-2}(u,v), \,(u,v)\in {\mathcal D}$ are $n-2$ mutually
orthogonal unit vectors, orthogonal to $\span \{l, l_u, l_v\}$.

\vskip 1mm

The minimal foliated semi-symmetric hypersurfaces in $\E^{n+1}$
are characterized analytically in  \cite {GanMil-CR}  by the
following

\vskip 2mm
\begin{thm}\label{T:basic}
Let $M^n$ be a hypersurface in $\E^{n+1}$ which is the envelope of
a two-parameter family of hyperplanes, determined by a unit
vector-valued function $l(u,v)$, represented by isothermal
parameters, and a scalar function $r(u,v)$. Then $M^n$ is minimal
if and only if \, $l(u, v)$ and $r(u,v)$ satisfy the equalities
$$\begin{array}{l}
\vspace{2mm}
 \Delta l(u,v) + 2E(u,v)\,l(u,v) = 0;\\
\vspace{2mm}
 \Delta r(u,v) + 2E(u,v)\,r(u,v) = 0.
\end{array}$$
\end{thm}

Hence, the minimal foliated semi-symmetric hypersurfaces in $\R^4$
are generated by the solutions of system (2.3), i.e. by the
minimal surfaces in $S^3$.

\vskip 1mm Now we shall construct examples of minimal foliated
semi-symmetric hypersurfaces in $\R^4$, which are generated by
some minimal surfaces in $S^3$.

\vskip 1mm The simplest example of a minimal surface in $S^3$ is
the sphere $S^2 = S^3 \bigcap \R^3$. We assume that $\R^3$ is the
subspace of $\R^4$ orthogonal to $e_4$, i.e. $\R^3 = \span\{e_1, e_2, e_3\}$. An isothermal parametrization of $S^2$ is given by
$$S^2: l(u,v) = \ds{\frac{1}{\cosh u} \left( \cos v;\, \sin v; \sinh u; \, 0 \right)}.$$
A direct computation shows that $E = \langle l_u, l_u \rangle =
\langle l_v, l_v \rangle = \ds{\frac{1}{\cosh ^2 u}}$, $F =
\langle l_u, l_v \rangle = 0$ and obviously $l(u,v)$ satisfies the
equality
$$ \Delta l(u,v) + \ds{\frac{2}{\cosh^2 u}}\, l(u,v) = 0.$$
The normal vector field $n(u,v)$ of $S^2$ is $n = e_4 = \left(0;\,
0;\, 0;\,1 \right)$. According to Theorem \ref{T:basic} the
corresponding differential equation for the scalar function
$r(u,v)$ is
$$ \Delta r(u,v) + \ds{\frac{2}{\cosh^2 u}}\, r(u,v) = 0. \leqno{(3.2)}$$
Every solution $r(u,v)$ of (3.2) together with the sphere $S^2: l =
l(u,v)$ generate a minimal foliated semi-symmetric hypersurface in
$\R^4$ according to formula (3.1).

One solution of (3.2) is
$$r(u,v) = \ds{\left(v+\frac{\pi}{2}\right) \tanh u}.$$
Now let us consider the  minimal foliated semi-symmetric
hypersurface $M^3$ generated by $l(u,v)$ and this solution
$r(u,v)$. Calculating $r_u$, $r_v$, $l_u$, $l_v$ and applying
formula (3.1), we obtain
$$M^3: X(u,v,w) = \left( -\sinh u\, \sin v;\, \sinh u\, \cos v;\,v+\frac{\pi}{2}; \, w \right).$$
After the following change of the parameters $u^1 = \sinh u; \,\,
t = \ds{v+\frac{\pi}{2}}$, we obtain the hypersurface
$$M^3: X(u^1,t,w) = u^1 (\cos t \,e_1 + \sin t \,e_2) + t e_3 + w e_4. \leqno{(3.3)}$$

The hypersurface $M^3$, whose radius-vector $X = X(u^1,t,w)$ is
determined by (3.3), is the so-called \textit{generalized
helicoidal ruled hypersurface}, obtained by  G. Aumann (see
\cite[Theorem 4]{A-min}). It is called also a \textit{first type
helicoid} in $\R^4$. It is a generalization of the  \textit{right
helicoid} in $\R^3$.

\vskip 2mm The next well-known example of a minimal surface in
$S^3$ is the Clifford torus
$$\mathcal{M}:  l(u,v) = (\cos u \,\cos v; \cos u\, \sin v; \sin u \,\cos v; \sin u\,\sin v).$$
The normal vector field $n(u,v)$ of $\mathcal{M}$ is
$$n(u,v) = (\sin u \,\sin v; -\sin u\, \cos v; -\cos u \,\sin v; \cos u\,\cos v).$$
Since $E = \langle l_u, l_u \rangle = \langle l_v, l_v \rangle
=1$, then according to Theorem \ref{T:basic} the corresponding
equation for the scalar function $r(u,v)$ is
$$ \Delta r(u,v) + 2\, r(u,v) = 0. \leqno{(3.4)}$$
If we take the trivial solution $r(u,v) = 0$ of (3.4), we get the
following minimal foliated semi-symmetric hypersurface
$$M^3: X(u,v,w) = w \,n(u,v) = w (\sin u \,\sin v; -\sin u\, \cos v; -\cos u \,\sin v; \cos u\,\cos v).$$
After the following change of the parameters $u^1 = -w \sin u;
\,\, u^2 = w \cos u; \,\, t = \ds{v+\frac{\pi}{2}}$, we obtain the
hypersurface
$$M^3: X(u^1,u^2,t) = u^1 (\cos t \,e_1 + \sin t \,e_2) + u^2 (\cos t \,e_3 + \sin t \,e_4). \leqno{(3.5)}$$

The hypersurface $M^3$, whose radius-vector $X = X(u^1, u^2,t)$ is
determined by (3.5), is the minimal ruled hypersurface obtained by
G. Aumann in  \cite[Theorem 1]{A-min}. This hypersurface is known
as a \textit{second type helicoid}. The first and the second type
helicoids are the only minimal ruled hypersurfaces in $\R^4$ (see
\cite{A-min}).

Thus we showed that the first type helicoid is generated by the
sphere $S^2$ in $S^3$, while the second type helicoid is generated
by the Clifford torus.

\vskip 1mm
 Our scheme of constructing minimal foliated
semi-symmetric hypersurfaces in $\R^4$ can be applied to each
minimal surface $M^2: l =l(u,v)$ in $S^3$ and each solution of the
corresponding differential equation for the scalar function
$r(u,v)$.

We shall illustrate how this construction can be applied to the
generalized torus of second type, given in Theorem \ref{T:Second
type}, in the special case when $t = 0$, $s = \ln \alpha,\,\,
\alpha>0, \, \alpha \neq 1$.
Since the calculations are too long and complicated, we give only a short sketch of the construction.
In this case the solution $l(u,v)$ is defined by
$$l(u,v) = \ds{f(u)\left[p(u)+ \frac{\alpha}{\alpha^2 + 1} \cos \beta v \left(\sqrt{\alpha} e_1 + \frac{1}{\sqrt{\alpha}} e_4\right) +
 \sqrt{\frac{\alpha}{\alpha^2 + 1}} \,\sin \beta v \,\,e_3\right]},$$
where
$f(u) = \ds{\sqrt{\frac{\alpha^2 \cos^2 h^{-1}(u) + \sin ^2 h^{-1}(u)}{\alpha}}}$,
$h^{-1}(u)$ is the inverse function of $h$ given on page \pageref{key}, $p(u)$ is the solution of system (2.42), and
$\beta = \ds{\sqrt{\frac{\alpha^2 + 1}{\alpha}}}$.
As a solution of the corresponding differential equation for $r(u,v)$ we choose
$r(u,v) = \ds{f(u)\sqrt{\frac{\alpha}{\alpha^2 + 1}}\,  \sin \beta v}$.
We calculate $l_u(u,v)$, $l_v(u,v)$, and using (2.33) and (2.40) we find the vector-valued function $n(u,v)$:
$$\begin{array}{ll}
\vspace{2mm}
n(u,v) = & \ds{\frac{1-\alpha^2}{\alpha(\alpha^2 + 1)} (e_1 -
\alpha e_4) + \frac{\alpha}{(\alpha^2 + 1) f(u)} \cos \beta
v\left(\sqrt{\alpha} e_1 + \frac{1}{\sqrt{\alpha}} e_4\right)}\\
\vspace{2mm}
& +  \ds{\sqrt{\frac{\alpha}{\alpha^2 + 1}}\frac{1}{f(u)} \,\sin \beta v \,\,e_3 - \frac{p(u)}{f(u)} -
\frac{1-\alpha^2}{\alpha} \int_0^u \frac{\sin 2h^{-1}(s)}{f^2(s)}\, p(s) ds}.
\end{array}$$

Applying formula (3.1) we obtain the following minimal foliated
semi-symmetric hypersurface $M^3$:
$$\begin{array}{ll}
\vspace{2mm}
X(u,v,w) = & \ds{\frac{1-\alpha^2}{2\sqrt{\alpha (\alpha^2 + 1)}} \, \frac{\sin 2h^{-1}(u)}{f(u)}\, \sin \beta v\, p'(u) }\\
\vspace{2mm} & \ds{+ \sqrt{\frac{\alpha}{\alpha^2+1}} \,\frac{\alpha^4 \cos^2 h^{-1}(u) + \sin^2 h^{-1}(u)}{\alpha^2 f^2(u)}\,  \sin \beta v\, p(u)}\\
\vspace{2mm}
& \ds{- \left(\sqrt{\frac{\alpha}{\alpha^2+1}}\right)^3 \frac{1}{2 f^2(u)} \, \sin 2 \beta v \left(\sqrt{\alpha} e_1 + \frac{1}{\sqrt{\alpha}} e_4\right)}\\
\vspace{2mm} & \ds{+  \left(1 - \frac{\alpha}{(\alpha^2+1) f^2(u)}\,\sin^2 \beta v \right) e_3 + w\, n(u,v)}.
\end{array} $$

Unfortunately, in this case we cannot write by means of elementary functions the hypersurface $M^3$ as in the previous examples, because
we cannot find  an explicit solution $p(u)$ of linear system (2.42).

\vskip 2mm \textbf{Acknowledgements:} The second author is
partially supported by "L. Karavelov" Civil Engineering Higher
School, Sofia, Bulgaria under Contract No 10/2009.


\begin{thebibliography}{99}

\bibitem{A-min}
Aumann G., {\it Die Minimalhyperregelfl\"{a}chen.} Mh. Math., {\bf
34} (1981), 293-304.

\bibitem{Enneper}
Enneper A., \textit{Die cyklischen Fl\"{a}chen}. Z. Math. Phys.
\textbf{14} (1869), 393-421.

\bibitem{GanMil-MathB}
Ganchev G. and Milousheva V.,  {\it A generation of foliated
semi-symmetric hypersurfaces in the four-dimensional Euclidean
space}. Mathematica Balkanica, {\bf 21} (2007) 1-2, 97-111.

\bibitem{GanMil-CR}
Ganchev G. and Milousheva V.,  {\it An analytic characterization
of the minimal and the bi-umbilical foliated semi-symmetric
hypersurfaces in Euclidean space}. C. R. Acad. Bulg. Sci., {\bf
60} (2007) 6, 601-606.

\bibitem{G-M-S}
Giaguinta M., Modica G. and Sou\v{c}ek J., \textit{Cartesian
currents in the calculus of variations, I and II}. Berlin,
Springer, 1998.


\bibitem{Hild}
Hildebrandt S. \textit{Nonlinear elliptic systems and harmonic
mappings}. Proc. Beijing Symp. on Differential Geometry and
Differential Equations, 1980.

\bibitem{Jost-1}
Jost J., \textit{Harmonic maps between surfaces (Lecture Notes in
Mathematics 1062)}, Springer, Berlin, 1984.

\bibitem{Jost-2}
Jost J., \textit{Two-dimensional geometric variational problems}.
Wiley-Interscience, Chichester, 1991.

\bibitem{Lawson-2}
Lawson H.B., Jr. {\it Local rigidity theorems for minimal
hypersurfaces}. Ann. of Math., {\bf 89} (1969), 187-197.

\bibitem{Lawson}
Lawson H.B., Jr. {\it Complete minimal surfaces in $S^3$}. Ann. of
Math., {\bf 92} (1970), 335-374.

\bibitem{Lopez}
L\'{o}pez R., \textit{Cyclic surfaces of constant Gauss
curvature}. Houston Math. J. \textbf{27} (2001), 799-805.

\bibitem{Nitsche}
Nitsche J., \textit{Cyclic surfaces of constant mean curvature}.
Nachr. Akad. Wiss. Gottingen Math. Phys. II \textbf{1} (1989),
1-5.

\bibitem{Riemann}
Riemann B., \textit{\"{U}ber die Fl\"{a}chen vom Kleinsten Inhalt
be gegebener Begrenzung}, Abh. K\"{o}nigl. Ges. Wissensch.
G\"{o}ttingen, Mathem. Cl. \textbf{13} (1868), 329-333.

\bibitem{Struwe}
Struwe M., \textit{Plateau's problem and the Calculus of
Variations}. Mathematical Notes  35. Princeton Univ. Press,
Princeton, 1988.

\bibitem{Szabo}
Szabo Z., {\it Structure theorems on Riemannian spaces satisfying
$R(X,Y) \cdot R=0$, I. The local version}.  J. Differ. Geom., {\bf
17} (1982), 531-582.

\end{thebibliography}
\end{document}